\documentclass[12pt]{article} 
\usepackage{fancyheadings}
\usepackage{geometry}
\usepackage{fullpage}
\usepackage{epsfig}
\usepackage[usenames]{color}
\usepackage{multicol}
\usepackage{appendix}
\newtheorem{theorem}{Theorem}[section]
\newtheorem{lemma}[theorem]{Lemma}

\newtheorem{algorithm}[theorem]{Algorithm}

\newcommand{\qed}{\nobreak \ifvmode \relax \else
      \ifdim\lastskip<1.5em \hskip-\lastskip
      \hskip1.5em plus0em minus0.5em \fi \nobreak
      \vrule height0.75em width0.5em depth0.25em\fi}

\def\cL{{\mathcal L}}
\def\la{\lambda}
\begin{document}
\newcommand{\mbf}[1]{\mbox{\boldmath $#1$}}

\thispagestyle{empty}
\setcounter{page}{0}
\begin{center}
\Large{\bf{\sf \textcolor{Mahogany}{LASSO-Patternsearch Algorithm with Application
to\\ Ophthalmology and Genomic Data}\footnote[1]{To appear, Statistics and Its Interface}}}
\end{center}
\vspace{.3in}
\begin{center}
{\sf \textcolor{Mahogany}{Weiliang Shi}\footnote[2]{Research supported in part by
NIH Grant EY09946, NSF Grants DMS-0505636, DMS-0604572 and ONR Grant
N0014-06-0095}\\
{\tt shiw@stat.wisc.edu}\\
Department of Statistics, University of Wisconsin \\
1300 University Avenue, Madison WI 53706 \\
\vspace{.2in}
{\sf \textcolor{Mahogany}{Grace Wahba}\footnote[3]{Corresponding Author.
Research supported in part by
NIH Grant EY09946, NSF Grants DMS-0505636, DMS-0604572 and ONR Grant
N0014-06-0095}\\
{\tt wahba@stat.wisc.edu}\\
Department of Statistics, Department of Computer Science and Department of
Biostatistics and Medical Informatics, University of Wisconsin\\
1300 University Avenue, Madison WI 53706\\
\vspace{.2in}
\textcolor{Mahogany}{Stephen Wright}\footnote[4]{Research supported
in part by NSF Grants SCI-0330538, DMS-0427689, CCF-0430504, CTS-0456694,
CNS-0540147 and DOE Grant DE-FG02-04ER25627}\\
{\tt swright@cs.wisc.edu}\\
Department of Computer Science, University of Wisconsin\\
1210 West Dayton Street, Madison WI 53706\\
\vspace{.2in}
\textcolor{Mahogany}{Kristine Lee}\footnote[5]{Research supported
in part by NIH grants EY06594 and EY015286}, 
\textcolor{Mahogany}{Ronald Klein}\footnote[6]{Research support in part by NIH grants 
EY06594, EY015286 and Research to Prevent Blindness Senior Investigator Awards},
\textcolor{Mahogany}{Barbara Klein}\footnotemark[6]\\
{\tt (klee, kleinr, kleinb)@epi.ophth.wisc.edu}\\
Department of Ophthalmology and Visual Science, University of Wisconsin\\
610 Walnut St., Madison WI\\ 
\vspace{0.3in}
}
}
\end{center}
\normalsize
\newpage
\normalsize
\newpage
\title{\textbf{LASSO-Patternsearch Algorithm with Application to Ophthalmology 
and Genomic Data }}
\author{\Large{Weiliang Shi, Grace Wahba, Stephen Wright, }\\ \Large{Kristine Lee, Ronald Klein and Barbara Klein}\\\Large{University of Wisconsin - Madison}}
\maketitle
\begin{abstract}
The LASSO-Patternsearch algorithm is proposed as a two-step method 
to identify clusters or patterns of 
multiple risk factors for outcomes of interest in  demographic 
and genomic studies. 
The predictor variables are dichotomous or can be coded 
as dichotomous. 
Many diseases are suspected of having multiple 
interacting risk factors acting in concert, and it is of much 
interest to uncover higher order interactions or risk patterns when they exist. 
The patterns considered here are those that arise naturally from the log linear 
expansion of the multivariate Bernoulli density. 
The method is designed for the case where there is a possibly
very large number of candidate patterns but it is believed that 
only a relatively small number are important. 
A LASSO is used to greatly reduce the number
of candidate patterns, 
using a novel computational algorithm 
that can handle an extremely large number of unknowns
simultaneously. Then the patterns surviving the 
LASSO are further pruned in the framework of (parametric) generalized 
linear models. A novel tuning procedure based on 
the GACV for Bernoulli outcomes, modified  to act as a model selector,
is used at both steps. 
We first  applied the method to myopia data from the 
population-based Beaver Dam Eye 
Study, exposing physiologically interesting interacting risk 
factors. 
In particular, we found that for an older cohort
the risk of myopic changes in refraction for 
smokers is reduced by taking vitamins while the risk for 
non-smokers is independent of the ``taking vitamins"  variable.
This is in agreement with the general result that smoking 
reduces the absorption of vitamins, and certain vitamins 
have been associated with eye health.
We then applied the method to 
data from a generative model of Rheumatoid
Arthritis based on  Problem 3 from the Genetic Analysis Workshop 15,
successfully demonstrating its potential to  
efficiently recover higher 
order patterns from attribute vectors of length typical of genomic
studies. 
\end{abstract}

\section{Introduction}
We consider the problem which occurs in demographic 
and genomic studies when there 
are a large number of risk factors that potentially interact in 
complicated ways to induce elevated risk. The goal is to search for important
patterns of  multiple risk factors among a very large number of candidate 
patterns, with results that are easily interpretable. In this work 
the LASSO-Patternsearch algorithm (LPS)  is proposed for this task. 
All variables are binary, or have been dichotomized before the analysis, 
at the risk of some loss of information;  this allows the study of much higher 
order interactions than would be possible with risk factors with more than
several possible values, or with continuous risk factors.
Thus LPS may, if desired,   be used as a preprocessor
to select clusters of variables that are later  analyzed in 
their pre-dichotomized form, see \cite{zhang:wahba:lin:voelker:2004}.
Along with demographic studies, a particularly promising application of LPS
is to the analysis of patterns 
or clusters of SNPs (Single Nucleotide Polymorphisms) or other
genetic variables that are associated with 
a particular phenotype, when the attribute vectors 
are very large and there exists a very large number 
of candidate patterns. LPS is designed specifically for the 
situation where the number of candidate patterns may be very large, but 
the solution, which may contain high order patterns, is believed to be 
sparse. LPS  is based on  the log linear parametrization 
of the multivariate
Bernoulli distribution \cite{whittaker:1990} to generate 
all possible patterns, if feasible, or at least a 
large subset of all possible patterns 
up to some maximum order.   LPS begins with a LASSO algorithm (penalized 
Bernoulli likelihood with an $l_1$ penalty),  used with a new 
tuning score, BGACV.  BGACV is a modified version of 
the GACV score {\cite{xiang:wahba:1996} to target 
variable selection, as opposed to Kullback-Liebler distance,
which is the GACV target.
A novel numerical algorithm is developed specifically for this 
step, which can handle an extremely large number of basis functions
(patterns) simultaneously. This is in particular contrast to most of the 
literature in the area, which uses greedy or sequential algorithms. 
The patterns surviving this process are then entered into a parametric 
linear logistic regression to obtain the final model, 
where further sparsity may be enforced via
a backward elimination process using the BGACV score as a stopping 
criterion. Properties of LPS will be  examined via simulation, 
and in demographic data by scrambling 
responses to establish false pattern 
generation rates.

There are many approaches that can model data with binary covariates and binary responses, 
see, for example CART \cite{breiman:friedman:olshen:stone:1984}, LOTUS \cite{chan:loh:2004}, 
Logic regression \cite{ruczinski:kooperberg:leblanc:2003} and 
Stepwise Penalized Logistic Regression (SPLR) \cite{park:hastie:2008}. 
Logic regression is an adaptive regression methodology that constructs predictors 
as Boolean combinations of binary covariates. It uses simulated 
annealing to search through the high dimensional covariate space and 
uses five-fold cross validation and randomization based hypothesis testing 
to choose the best model size. SPLR is a variant of logistic 
regression with $l_2$ penalty to fit interaction models. It uses a forward 
stepwise procedure to search through the high dimensional covariate space. 
The model size is chosen by an  AIC- or BIC-like score and the 
smoothing parameter is chosen by 5-fold cross validation. 
For Gaussian data the LASSO was proposed in \cite{tibshirani:1996} as a variant of 
linear least squares ridge regression with many predictor variables. 
As proposed there, the LASSO  minimized the residual sum of squares 
subject to a constraint that the sum of absolute values of the coefficients 
of the basis functions be 
less than some constant, say $t$. This is equivalent to minimizing 
the residual sum of squares plus a penalty which is some multiple 
$\lambda$ (depending on $t$)  of the sum of absolute values ($l_1$ penalty). 
It was demonstrated there that this approach tended to set many of 
the coefficients to zero, resulting in a sparse model, 
a property not generally obtaining with quadratic penalties. 
A similar idea was exploited in \cite{chen:donoho:saunders:1998} 
to select a good subset of an overcomplete set of nonorthogonal 
wavelet basis functions. The asymptotic behavior of LASSO type estimators was 
studied in \cite{knight:fu:2000}, and \cite{osborne:presnell:turlach:2000} 
discussed computational procedures in the Gaussian context. 
More recently \cite{efron:hastie:johnstone:tibshirani:2004} discussed 
variants of the LASSO and methods for computing the LASSO for a 
continuous range of values of $\lambda$ in the Gaussian case. 
Variable selection properties of the LASSO were examined in 
\cite{leng:lin:wahba:2006} in some special cases, and many 
applications can be found on the web. In the context of nonparametric 
ANOVA decompositions \cite{zhang:wahba:lin:voelker:2004} used 
an overcomplete set of basis functions obtained from a Smoothing Spline 
ANOVA model, and used $\ell_1$ penalties on the coefficients of main 
effects and low order interaction terms, in 
the spirit of \cite{chen:donoho:saunders:1998}. 
The present paper uses some ideas from \cite{zhang:wahba:lin:voelker:2004},
although the basis function set here is quite different.
Other work
has implemented  
$\ell_1$ penalties along with quadratic (reproducing kernel square norm) 
penalties to take advantage of the properties 
of both kinds of penalties, 
see for example 
\cite{gunn:kandola:2002} \cite{lee:kim:lee:koo:2006} \cite{zhang:lin:2006}\cite{zou:hastie:2005b}.

The rest of the article is organized as follows. In Section 2 we describe 
the first (LASSO)  
step of the LPS
including choosing the smoothing parameter by the 
B-type Generalized Approximate Cross Validation (BGACV), ``B"  standing for the 
prior belief that the solution is sparse, analogous to BIC.  An 
efficient algorithm for the LASSO step is presented here. 
Section 3  describes the second step of 
the LASSO-Patternsearch algorithm, utilizing a parametric 
logistic regression, again tuned by BGACV.
Section 4 presents three simulation examples, designed to demonstrate
the properties of LPS as well as comparing LPS to Logic regression 
and SPLR. Favorable properties of LPS are exhibited in models 
with high order patterns and correlated attributes. 
Section 5 applies the method to myopic changes in refraction in
an older cohort from the Beaver Dam Eye Study 
\cite{klein:klein:linton:demets:1991}, where some interesting
risk patterns including one involving smoking and vitamins are found.
Section 6 
applies the method to data from a generative
model of Rheumatoid Arthritis Single Nucleotide Polymorphisms
adapted from the Genetics Analysis Workshop 15 \cite{maccluer:2007}, 
which examines
the ability of the algorithm to recover third order patterns
from extremely large attribute vectors.
Section 7 
notes some generalizations, and, finally, Section 8  gives a summary and 
conclusions. 
Appendix A derives the BGACV score;
Appendix B gives details of the specially designed code for the LASSO which is 
capable of handling a very large number of patterns simultaneously; 
Appendix C shows the detailed results of Simulation Example 3. 
When all of the variables are coded as 1 in the risky direction, 
the model will be sparsest among equivalent models. 
Appendix D gives a lemma describing what happens when 
some of the variables are coded with the opposite direction as 1.

\section{The LASSO-Patternsearch Algorithm}
\subsection{The LASSO-Patternsearch Algorithm - Step 1}
Considering $n$ subjects, for which $p$ variables are observed, we first reduce  continuous variables to ``high" or ``low" in order to be able to examine {\it very many} variables and {\it their interactions}  simultaneously. We will assume that for all or most of the the $p$ variables, we know in which direction they are likely to affect the outcome or outcomes of interest, if at all. For some variables, for example smoking,  it is clear for most endpoints in which direction the smoking variable is likely to be  ``bad" if it has any effect, and this is true of many but not all variables. For some continuous variables, for example systolic blood pressure, higher is generally ``worse", but extremely low can also be ``bad". For continuous variables, we need to initially assume the location of a cut point on  one side of which the variable is believed to be ``risky" (``high")  and the other side ``not risky" (``low"). For systolic blood pressure that might, for example, be 140 mmHg.  For an economic variable,
that might be something related to the current standard for poverty level. If the ``risky" direction is known for most variables the results will be readily interpretable. Each subject thus has an 
attribute vector of $p$ zeroes and ones, describing whether each of their $p$ attributes is on one side or the other of the cutoff point. The LASSO-Patternsearch approach described below is able to deal with high order interactions and very large $p$. 
The data is $\{y_i, x(i), i = 1,  \cdots , n\}$, where 
$y_i \in \{0,1\}$ codes the response, 
$x(i) = (x_1(i), x_2(i), \cdots , x_p(i)) $ is the attribute vector for the $i$th subject, $x_j(i) \in \{0,1\}$. Define the basis functions $B_{j_1 j_2..j_r}(x) = \prod_{\ell= 1}^r x_{j_{\ell}}$, that is, $B_{j_1 j_2..j_r}(x)  = 1$ if $x_{j_1}, ..., x_{j_r}$ are all $1$'s and $0$ otherwise. We will call $B_{j_1 j_2..j_r}(x)$ an $r$th order pattern. 
Let $q$ be the highest order we consider.
Then there will be $N_B = \sum_{\nu =0}^{q} {{p}\choose{\nu}}$ 
patterns. If $q = p$, we have a complete set of $N_B = 2^p  $ such 
patterns (including the constant function $\mu$), spanning all possible patterns. 
If $q = 1$ only first order patterns
(henceforth called ``main effects") are considered, if 
$q = 2$ main effects and second order 
patterns are considered, and so forth.  
Letting $p(x) = \mbox{Prob}[y= 1|x]$ 
and the logit (log odds ratio) be $f(x) = \log [p(x)/(1-p(x))]$, 
we estimate $f$ by minimizing 
\begin{equation}\label{penlik}
I_{\lambda}(y,f) =  \cL (y,f) +\lambda J(f)
\end{equation}
where $\cL (y,f)$ is $\frac{1}{n}$ times the negative log likelihood:
\begin{equation}\label{defL}
\cL(y,f) =\frac{1}{n}  \sum_{i=1}^n [ -y_if(x(i)) + \log (1+e^{f(x(i))})]
\end{equation}
with 
\begin{equation}\label{f.def}
f(x) = \mu +  \sum_{\ell=1}^{N_B-1} c_{\ell}B_{\ell}(x),
\end{equation}
where we are relabeling the $N_B-1$ (non-constant) patterns from $1$ to $N_B-1$, and 
\begin{equation}\label{J}
J(f) =  \sum_{\ell = 1}^{N_B-1} |c_{\ell}|.
\end{equation}
If all possible patterns are included in (\ref{f.def})
then $f$ there is the most  general form of the log odds
ratio for $y$ given $x$ obtainable from the log linear parametrization
of the multivariate Bernoulli distribution given in \cite{whittaker:1990}.
In Step 1 of the LASSO-Patternsearch we minimize (\ref{penlik}) 
using  the  BGACV score to choose $\lambda$. 
The next section describes the BGACV 
score and the kinds of results it can be expected to produce.
\subsection{B-type Generalized Approximate Cross Validation (BGACV)}\label{btype}
The tuning parameter $\lambda$ in ($\ref{penlik}$) balances the trade-off 
between data fitting and the sparsity of the model. 
The bigger $\lambda$ is, the sparser the model. 
The choice of $\lambda$ is generally a crucial part of 
penalized likelihood methods and machine learning techniques 
like the Support Vector Machine. For smoothing spline models 
with Gaussian data, \cite{wahba:wold:1975} proposed ordinary 
leave-out-one cross 
validation (OCV)  . Generalized Cross Validation (GCV), derived 
from OCV,  was proposed in 
\cite{craven:wahba:1979}\cite{golub:heath:wahba:1979}, 
and theoretical properties were obtained in \cite{li:1986}
and elsewhere.  
For smoothing spline models with Bernoulli data and quadratic 
penalty functionals, \cite{xiang:wahba:1996} derived the
Generalized Approximate Cross Validation (GACV) 
from an OCV estimate following the method used to obtain GCV. 
In \cite{zhang:wahba:lin:voelker:2004} GACV  was extended to the 
case of Bernoulli data with continuous covariates 
and $l_1$ penalties. 

The derivation of the GACV  begins with a leaving-out-one 
likelihood to minimize an estimate of the comparative 
Kullback-Leibler distance (CKL) between the true and 
estimated model distributions. The ordinary leave-out-one 
cross validation score for CKL is 
\begin{equation}
CV(\lambda) = \frac{1}{n} \sum_{i=1}^n 
[-y_if_{\lambda}^{[-i]} (x(i)) + \log (1+e^{f_{\lambda}(x(i))})],
\end{equation}
where $f_{\la}$ is the minimizer of the objective function 
(\ref{penlik}), and $f_{\lambda}^{[-i]}$ is the minimizer of 
(\ref{penlik}) with the $i$th data point left out. 
Through a series of approximations and an averaging step 
as described in Appendix A, we obtain the GACV score 
appropriate to  the present context. It is a simple to 
compute special case
of the GACV score in \cite{zhang:wahba:lin:voelker:2004}:

\begin{equation}\label{gacv}
GACV(\lambda) = \frac{1}{n}\sum_{i=1}^n[-y_if_{\lambda i} 
+ \mbox{log}(1+e^{f_{\lambda i}})] + 
\frac{1}{n} {\mbox{tr}}H
\frac{\sum_{i=1}^ny_i(y_i-p_{\lambda i})}{(n - N_{B_0}) },
\end{equation}
here $H = B^*({B^*}'WB^*)^{-1}{B^*}'$, where W is the $n \times n$ 
diagonal matrix with $ii$th element 
the estimated 
variance at $x(i)$ 
($p_{i\lambda}(1-p_{i\lambda}))$
and $B^*$ is the $n \times N_{B_0}$ 
design matrix for the $N_{B_0}$ non-zero $c_{\ell}$ in the model.
The quantity 
${\mbox{tr}}H
\frac{\sum_{i=1}^ny_i(y_i-p_{\lambda i})}{(n - N_{B_0}) }$
plays the role of degrees of freedom here.
As is clear from the preceding discussion, the GACV is a criterion
whose target is the minimization of the (comparative) Kullback-Liebler
distance from the estimate to the unknown ``true" model. By analogy with the 
Gaussian case (where the predictive mean squared error and comparative
Kullback-Liebler distance coincide), it is known that
optimizing for predictive mean square error and optimizing for 
model selection {\it when the true model is sparse} are 
not in general the same thing. This is discussed in various 
places, for example, see  \cite{george:2000} 
which discusses the relation 
between AIC and BIC, AIC being a predictive criterion and BIC, 
which generally results in a sparser model, being a model 
selection criterion, with desirable properties when the 
``true" model is of fixed (low) dimension as the sample size 
gets large. See also
\cite{haughton:1988}, \cite{leng:lin:wahba:2006} and 
particularly our remarks at the end of Appendix \ref{BGACV}. 
In the  AIC to BIC transformation, if $\gamma$ is the degrees of 
freedom for the model, then BIC replaces $\gamma$ with 
$(\frac{1}{2}\log n)\gamma$. By analogy we obtain  
a model selection criterion, BGACV,  from GACV as follows.  
Letting $\gamma$ be the quantity playing the role of degrees
of freedom for signal in the Bernoulli-$l_1$ penalty case,  
\begin{equation}\label{q}
\gamma = {\mbox{tr}}H
\frac{\sum_{i=1}^ny_i(y_i-p_{\lambda i})}{(n - N_{B_0}) },
\end{equation}
$\gamma$ is replaced by $(\frac{1}{2}\log n)\gamma$ to 
obtain
\begin{equation}\label{bgacv}
BGACV(\lambda) = \frac{1}{n}\sum_{i=1}^n[-y_if_{\lambda i} 
+ \mbox{log}(1+e^{f_{\lambda i}})] 
+\frac{1}{n} 
\frac{\log n}{2}
\mbox{tr}H
\frac{\sum_{i=1}^ny_i(y_i-p_{\lambda i})}{(n - N_{B_0}) }.
\end{equation}

We illustrate the difference of empirical performances between GACV and BGACV 
on a ``true" model with a small number of strong patterns.
Let $(X_1^*,X_4^*)^T$, $(X_2^*,X_5^*)^T$ and $(X_3^*,X_6^*)^T$ 
be independently distributed from a bivariate normal 
distribution with mean 0, 
variance 1 and covariance 0.7. $X_i = 1$ if $X_i^*>0$ and 0 otherwise, 
$i = 1,2,\cdots, 6$. $X_7$ is independent of the others and takes 
two values $\{1,0\}$, each with a probability of 0.5. 
$X = (X_1, \cdots, X_7)$. 
The sample size $n = 800$. 
Three patterns that consist of six variables are 
important, and $X_7$ is noise. 
The true logit is
\begin{equation}
f(x) = -2 + 1.5B_1(x) + 1.5B_{23}(x)+2B_{456}(x).\nonumber
\end{equation}
This problem is very small so we chose the maximum order $q = p = 7$, so 
127 patterns plus a constant term are entered in the trial model. 
We ran the simulation 100 times and the result is shown in Table \ref{sim1tab}. 
Both GACV and BGACV select the three important patterns perfectly, but 
GACV selects more noise patterns than BGACV. Note that a total of $2^7 - 4$ 
noise patterns have been considered in each run. The maximum possible 
number in the last column is $100\times (2^7-4) = 12,400$. 
Neither GACV or BGACV is doing a bad job but we will discuss a method to 
further reduce the number of  noise patterns in Section 3. 
Figure \ref{cvplot} shows the scores of these two criteria in the 
first data set. BGACV selects a bigger smoothing parameter than GACV does. 
These scores are not continuous at the the point where a parameter 
becomes zero so we see jumps in the plots. 

\begin{table}[htbp]
  \caption{The results of Simulation Example 1. The second through 
fourth columns are the appearance frequencies of the three important 
patterns in the 100 runs. The last column is the total appearance frequency 
of all other patterns. The second and third row compare GACV and BGACV 
within the first step of LPS. The fourth through sixth rows, 
which compare the full LPS with Logic regression and SPLR will be
discussed in Section \ref{simex1}.} 
  \begin{center}
    \begin{tabular}{c|cccc}
      \hline
      Method&$B_1$&$B_{23}$&$B_{456}$&other\\\hline
      $~~~$GACV$~~~$&$~~~~100~~~~$&$~~~~100~~~~$&$~~~~100~~~~$&$~~~~749~~~~$\\
      BGACV&100&100&100&568\\\hline
      LPS&97&96&98&34\\
      Logic&100&94&94&64\\
      SPLR&100&90&55&426\\\hline
    \end{tabular}
  \end{center}
  \label{sim1tab}
\end{table}

\begin{figure}
  \begin{center}
    \includegraphics[scale = 0.6]{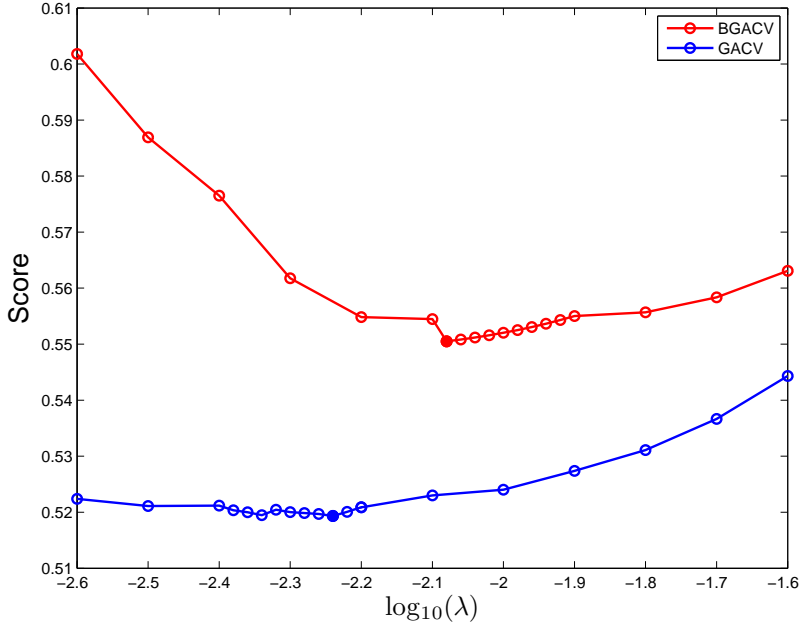}\\
  \end{center}
  \caption{Comparison of BGACV and GACV in the first data set of Simulation Example 1. The solid dots are the minima. BGACV selects a bigger $\lambda$ than GACV does.}
\label{cvplot}
\end{figure}
\subsection{Computation}
From a mathematical point of view, this optimization problem (\ref{penlik}) is 
the same as the likelihood basis pursuit (LBP) algorithm in \cite{zhang:wahba:lin:voelker:2004}, 
but with different basis functions. The solution can easily be computed via a general constrained nonlinear minimization code such as MATLAB's {\tt fmincon} on a desktop, for a range of values of $\lambda$, provided $n$ and
$N_B$ are not too large. However, for extremely large data sets 
with  more than a few attributes $p$ (and therefore a large number $N_B$ of possible basis functions), the problem becomes much more difficult to solve computationally with general optimization software, and algorithms that exploit the structure of the problem are needed. We design an algorithm that uses gradient information for the likelihood term in (\ref{penlik}) to find an estimate of the correct active set (that is, the set of components $c_\ell$ that are zero at the minimizer).  When there are not too many nonzero parameters, the algorithm also attempts a Newton-like enhancement to the search 
direction, making use of the fact that first and second partial derivatives of the function in (\ref{penlik}) with respect to the coefficients $c_{\ell}$ are easy to compute analytically once the function has been evaluated at these values of $c_{\ell}$. It is not economical to compute the full Hessian (the matrix of second partial derivatives), so the algorithm 
computes only the second derivatives of the log likelihood function $\cL$ with respect 
to those coefficients $c_{\ell}$ that appear to be nonzero at the solution. 
For the problems that the LASSO-Patternsearch is designed 
to solve, just a small fraction of these $N_B$ coefficients are nonzero at the solution. 
This approach is similar to the two-metric gradient projection approach for 
bound-constrained minimization, but avoids duplication of variables and 
allows certain other economies in the implementation.

The algorithm is particularly well suited to solving the problem (\ref{penlik}) for a 
number of different values of $\lambda$ in succession; the solution for 
one value of $\lambda$ provides an excellent starting point for 
the minimization with a nearby value of $\lambda$. Further details of 
this approach can be found in Appendix \ref{algorithm}.
\section{The LASSO-Patternsearch Algorithm - Step 2}
In Step 2 of LASSO-Patternsearch algorithm, the $N_{B_0}$ patterns surviving Step 1 are entered into a linear logistic regression model 
using {\tt glmfit} in MATLAB 
and 
pattern selection is then carried out by the backward elimination method. We take out one of the $N_{B_0}$ patterns at a time, fit the model with the 
remaining patterns and compute the tuning score. 
The pattern that gives the best 
tuning score to the model after being taken out is removed from the model. This process continues until there are no patterns in the model. A final model is 
chosen from the pattern set with the best tuning score. 
Note that all-subset selection is not being done, 
since this will introduce an overly large number of degrees of 
freedom into the process.

If copious data is available, then a tuning set can be 
used to create the tuning score, but  
this is frequently not the case.
Inspired by the tuning method in the LASSO step, we propose the 
BGACV score for the parametric logistic regression.  The likelihood function is smooth with respect to the parameters so the robust assumption that appears in Appendix \ref{BGACV} is not needed.  Other than that, the derivation of the BGACV score for parametric logistic regression follows that in Appendix \ref{BGACV}.  Let $s$ be the current subset of patterns under consideration and $B_s$ be the design matrix. The BGACV score for logistic regression is the same as (\ref{bgacv}) with $H = B_s(B_s'WB_s)^{-1}B_s'$. 

\begin{equation}\label{bgacvp}
BGACV(s) = \frac{1}{n}\sum_{i=1}^n[-y_if_{s i} 
+ \mbox{log}(1+e^{f_{s i}})] + 
\frac{1}{n}\frac{\mbox{log}n}{2} \mbox{tr}H
\frac{\sum_{i=1}^ny_i(y_i-p_{s i})}{(n - N_{B_0}) },
\end{equation}
where $f_{s i}$ is the estimated log odds ratio for observation 
$i$ and $p_{si}$ is the corresponding probability. The BGACV scores are 
computed for each model that is considered in the backward elimination procedure, and the model with the smallest BGACV score is taken as the final model.

The following is a summary of the  LASSO-Patternsearch algorithm:

\begin{enumerate}
  \item Solve ($\ref{penlik}$) and choose $\lambda$ by BGACV. 
Keep the patterns with nonzero coefficients. 
 \item Put the patterns with nonzero coefficients from Step 1 into a logistic regression model and select models by the backward elimination method with the selection criterion being BGACV.
\end{enumerate}
For simulated data, the results can be compared with the 
simulated model. For observational data, a selective 
examination of the data will be used to validate the results.
Other logistic regression codes, e. g. from R or SAS can be 
used instead of {\tt glmfit} here.
\section{Simulation Studies}
In this section we study the empirical performance of the LPS through three simulated examples. The first example continues with simulated data in Section \ref{btype}.  
There are three pairs of correlated variables and one independent variable. 
Three patterns are related to the response. The second example has only one high order pattern. The correlation within variables in the pattern is high and the correlation between variables in the pattern and other variables varies. The last example studies the performance of our method under various correlation settings. 
We compare LPS with two other methods, 
Logic regression \cite{ruczinski:kooperberg:leblanc:2003} and 
Stepwise Penalized Logistic Regression (SPLR) \cite{park:hastie:2008}. 
We use the R package {\tt LogicReg} to run Logic regression 
and the R package {\tt stepPlr} to run SPLR. The number of 
trees and number of leaves in Logic regression are selected by 
5-fold cross validation. The smoothing parameter 
in SPLR is also selected by 5-fold cross validation, and 
then the model size is selected by a BIC-like score based
on an approximation to a  degrees of freedom reproduced in the 
Comments section of Appendix \ref{BGACV}. 

\subsection{Simulation Example 1}\label{simex1}
In this example we have 7 variables and the sample size is 800. The true logit is $f(x) = -2 + 1.5B_1(x) + 1.5B_{23}(x) + 2B_{456}(x)$. The distribution of the covariates was described in 
Section \ref{btype}.  We simulated 100 data sets according to this model and ran all three methods on these data sets. The results are shown in the last three rows of Table \ref{sim1tab}. 

Let's compare LPS with the LASSO step (third row in Table \ref{sim1tab}) first. LPS misses all three patterns a few times. However, these numbers are still very close to 100 and more importantly, LPS significantly reduced the number of noise patterns, from over 500 to 34. Here we see why a second step is needed after the LASSO step. Now let's look at LPS compared with the other two methods. Logic regression picks the first term perfectly but it doesn't do as well as LPS on the remaining two patterns. It also selects more noise patterns than LPS. SPLR does worse, especially on the last pattern. 
It is not surprising because this example is designed to be 
difficult for SPLR, which is a sequential method. 
In order for $B_{456}$ to be in the model, at least one main effect of $X_4$, $X_5$ and $X_6$ should enter 
the model first, say $X_4$. And then a second order pattern should also enter 
before $B_{456}$. It could be $B_{45}$ or $B_{46}$. However, none of these lower order patterns 
are in the true model. This makes it very hard for SPLR to consider $B_{456}$, and the 
fact that variables in the higher order pattern are correlated with variables in the lower order patterns makes it even harder. We also notice that SPLR 
selects many more noise patterns than LPS and Logic regression. 
Because of the way it allows high order patterns to enter the model, the 
smallest SPLR model has 6 terms, $B_1$, one of $B_2$ and $B_3$, $B_{23}$, one main 
effect and one second order pattern in $B_{456}$, and $B_{456}$. Conditioning on the appearance 
frequencies of the important patterns, SPLR has to select at least $90 + 2\times 55 = 200$ 
noise patterns. The difference, $426-200 = 226$ is still much bigger than  34 selected  by LPS.

\subsection{Simulation Example 2}
We focus our attention on a high order pattern in this example. Let $X_1^*$ through $X_4^*$ be 
generated from a normal distribution with mean 1 and variance 1. 
The correlation between any of these two is 0.7. $X_i=1$ if $X_i^*>0$ and 0 otherwise for $i = 1,\cdots,4$. $X_{i+4} = X_i$ with probability $\rho$ and $X_{i+4}$ will be generated from Bernoulli(0.84) otherwise for $i = 1,\cdots,4$. $\rho$ 
takes values $0, 0.2, 0.5$ and $0.7$. 
$X = (X_1, \cdots, X_8)$.
Note that $P(X_1=1) = 0.84$ in our simulation. 
The sample size is 2000 and the true logit  is 
$f(x) = -2 + 2B_{1234}(x)$. We consider all possible patterns, so $q = p = 8$. We also ran this example 100 times and the results are shown in Table \ref{sim2tab}. 

LPS does a very good job and it is very robust against 
increasing  $\rho$, which governs the correlation between 
variables in the model and other variables. 
We selected the high order pattern almost 
perfectly and kept the noise patterns below 10 in all four settings. Logic regression selects the important pattern from 70 to 80 times and noise patterns over 130 times. There is a mild trend that it does worse as the correlation goes up, but the last one is an exception. SPLR is robust against the correlation but it doesn't do very well. It selects the important pattern from 50 to 60 times and noise patterns over 500 times. From this example we can see that LPS is extremely powerful in selecting high order patterns. 

\begin{table}
\caption{The results of simulation example 2. The numerators are the appearance frequencies of $B_{1234}$ and the denominators are the appearance frequencies of all noise patterns.}
  \begin{center}
    \begin{tabular}{c|cccc}
      \hline
      $\rho$&0&0.2&0.5&0.7\\\hline
      LPS&98/10&100/9&97/8&98/5\\
      $~~$Logic$~~$&$~~~82/132~~~$&$~~~73/162~~~$&$~~~72/237~~~$&$~~~74/157~~~$\\
      SPLR&53/602&60/576&57/581&58/552\\\hline
    \end{tabular}
  \end{center}
  \label{sim2tab}
\end{table}

\subsection{Simulation Example 3}
The previous two examples have a small number of variables so we considered patterns of all orders. To demonstrate the power of our algorithm, we add in more noise variables in this example. The setting is  
similar to Example 2. Let $X_1^*$ through $X_4^*$ be generated from a 
normal distribution with mean 1 and variance 1. The correlation between 
any of these two is $\rho_1$ and $\rho_1$ takes values 
in 0,0.2, 0.5 and 0.7. $X_i=1$ if $X_i^*>0$ and 0 
otherwise, $i = 1,2,3,4$. $X_{i+4} = X_i$ with 
probability $\rho_2$ and $X_{i+4}$ will be generated from 
Bernoulli(0.84) otherwise for $i = 1,2,3,4$. $\rho_2$ takes values $0, 0.2, 0.5$ and $0.7$ 
also.  $X_9$ through $X_{20}$ are generated from Bernoulli(0.5) independently. 
$X = (X_1, \cdots, X_{20}).$
The sample size $n = 2000$ and the true logit is
\[f(x) = -2+2B_9(x) +2B_{67}(x) +2B_{1234}(x).\]
Unlike the previous two examples, we consider patterns only  up 
to the order of 4 because of the large number of variables. That gives us a total of ${{20}\choose{0}} + {{20}\choose{1}} +\cdots+{{20}\choose{4}} = 6196$ basis functions. 

\begin{figure}
  \begin{center}
    \includegraphics[scale = 0.8]{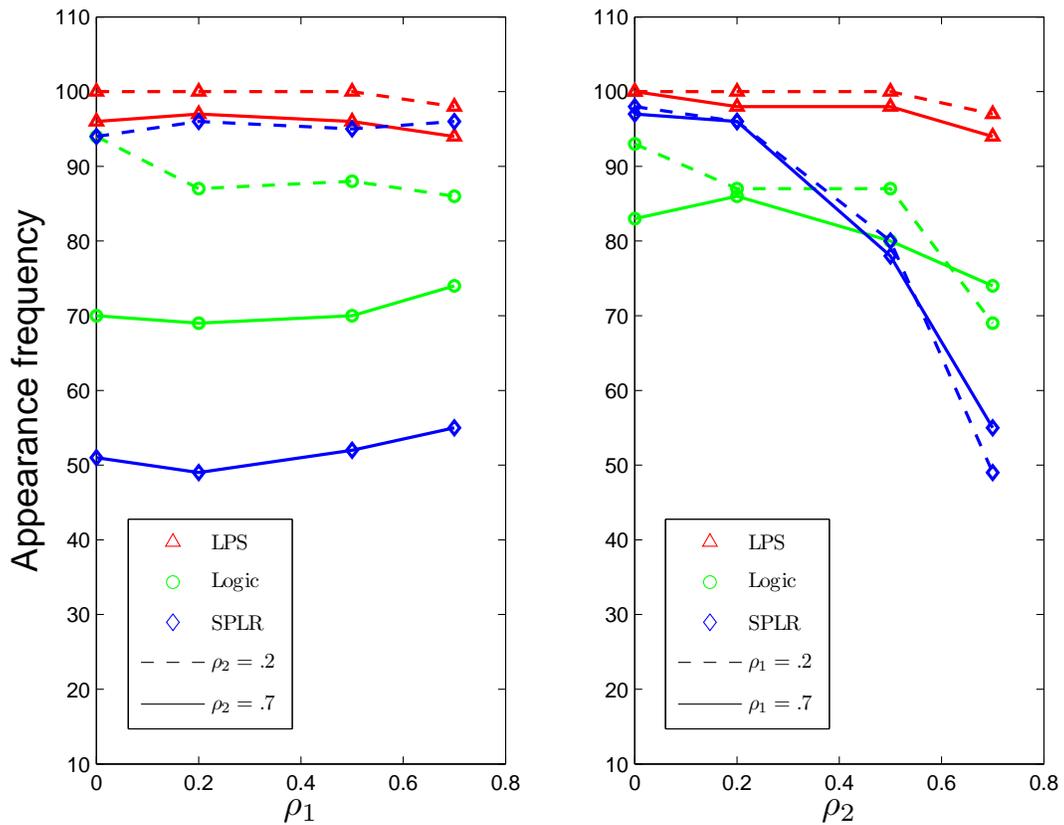}\\
  \end{center}
  \caption{Appearance frequency of the high order pattern $B_{1234}$ in simulation example 3. In the left panel, the x-axis is $\rho_1$. $\rho_2$ is 0.2 for the dashed line and 0.7 for the solid line. In the right panel, the x-axis is $\rho_2$. $\rho_1$ is 0.2 for the dashed line and 0.7 for the solid line. The red triangles represent LPS, the blue diamonds represent 
Logic regression \cite{ruczinski:kooperberg:leblanc:2003}
and the green circles represent SPLR \cite{park:hastie:2008}.}
\label{sim3plot}
\end{figure}

Figure \ref{sim3plot} shows the appearance frequencies of the high order pattern $B_{1234}$. From the left plot we see that LPS dominates the other two methods. All methods are actually very robust against $\rho_1$, the correlation within the high order pattern. There is a huge gap between the two blue lines, which means SPLR is very sensitive to $\rho_2$ which governs the correlation between variables 
in the high order pattern and others. This is confirmed by the right plot, where the blue lines decrease sharply as $\rho_2$ increases. We see similar but milder behavior in Logic regression. This is quite natural because the problem becomes harder as the the noise variables become more correlated with the important variables. However, LPS handles this issue quite well, at least in the current setting. We see a small decrease in LPS as $\rho_2$ goes up but those numbers are still very close to 100. The performance of these methods on the second order pattern $B_{67}$ is generally similar but the trend is less obvious as a lower order pattern is easier for most methods. The main effect $B_9$ is selected almost perfectly by every method in all settings. More detailed results are presented in Table \ref{sim3tab} in Appendix \ref{appsim3}. 

\section{The Beaver Dam Eye Study}
The Beaver Dam Eye Study (BDES) is an ongoing population-based study 
of age-related ocular disorders including cataract, age-related 
macular degeneration, visual impairment and refractive errors.
Between 1987 and 1988, a private census identified 5924 people aged 43 through 84 years in Beaver Dam, WI. 
4926 of these people participated the baseline exam (BD I) between 
1988 and 1990. Five (BD II), ten (BD III) and 
fifteen (BD IV) year follow-up data have been collected and 
there have been several hundred publications on this data. A detailed description of the study can be found in \cite{klein:klein:linton:demets:1991}. 

Myopia, or nearsightedness, is one of the most prevalent 
world-wide eye conditions.  
Myopia occurs when the eyeball is slightly longer than usual from front to back 
for a given level of refractive power of the cornea and lens and people with 
myopia find it hard to see objects at a distance
without a corrective lens. Approximately one-third of the 
population experience this eye problem and in some countries like Singapore, 
more than 70$\%$ of the population have myopia upon completing college \cite{seet:wong:tan:saw:2001}.
It is believed that myopia is related to various environmental 
risk factors as well as genetic factors. 
Refraction is the continuous measure from which myopia is 
defined.  Understanding how refraction changes over time can 
provide further insight into when myopia may develop.
Five and ten-year changes of refraction for the BDES population 
were summarized in \cite{lee:klein:klein:1999}\cite{lee:klein:klein:wong:2002}. 
We will study  five-year 
myopic changes in refraction (hereinafter called ``myopic change") 
in an older cohort 
aged 60 through 69 years. We focus on a small age group since 
the change of refraction differs for different age groups. 

Based on \cite{lee:klein:klein:wong:2002} 
and some preliminary analysis we 
carried out on this data, we choose seven risk 
factors: {\it sex, inc, jomyop, catct, pky, asa} and {\it vtm}
(sex, income, juvenile myopia, nuclear cataract, packyear, 
aspirin and vitamins). 
Descriptions and  binary cut points are presented in Table \ref{mypvar}. 
For most of these variables, we know which direction is bad. 
For example, male gender is a risk factor for most diseases and 
smoking is never good. The binary cut points are somewhat subjective here. Regarding {\it pky}, a pack a day for 30 years, for example, is a fairly substantial smoking history. {\it catct} has five levels of severity and we cut it at the third level. Aspirin $(asa)$ and vitamin supplements $(vtm)$ are 
commonly taken to maintain good health so we treat not taking 
them as risk factors. Juvenile myopia {\it jomyop} is assessed 
from self-reported age at which the person first started wearing 
glasses for distance. For the purposes of this study we have defined
myopic change  as a change in refraction of 
more than -0.75 diopters from baseline exam to the five year followup;
accordingly $y$ is assigned 1 if this change occurred
and 0 otherwise. 
There are 1374 participants in this age group at the baseline examination, 
of which 952 have measurements of refraction at the baseline and the five-year follow-up. Among the 952 people, 76 have missing values in the covariates. We assume that the missing values are missing at random for both response and covariates, although this assumption is not necessarily valid. However the examination of the missingness and possible imputation are beyond the scope of this study. Our final data consists of 876 subjects without any missing values in the seven risk factors. 

\begin{table}[htbp]
\caption{The variables in the myopic change example. The fourth column shows which direction is risky. }
  \begin{center}
    \begin{tabular}{c|lll}
      \hline
      code&variable&unit&higher risk\\\hline
      {\it sex}&sex&&Male\\
      {\it inc}&income& $\$$1000&$<$ 30\\
      {\it jomyop}&juvenile myopia&age first wore glasses for distance &yes before age 21\\
      {\it catct}&nuclear cataract&severity 1-5&4-5\\
      {\it pky}&packyear &pack per day$\times$years smoked&$>$30\\
      {\it asa}&aspirin&taking/not taking&not taking\\
      {\it vtm}&vitamins&taking/not taking&not taking\\\hline
    \end{tabular}
  \end{center}
  \label{mypvar}
\end{table}
As the data set is small, we consider all possible patterns ($q = 7$). 
The first step of the 
LASSO-Patternsearch algorithm selected 8 patterns, 
given in Table \ref{mypstep1pattern}.

\begin{table}[htbp]
\caption{Eight patterns selected at Step 1 in the myopic change data}
  \begin{center}
    \begin{tabular}{lll|llll}
      \hline
      &Pattern&Estimate $~~~$& $~~~$&& Pattern &Estimate\\\hline
      1&{\it constant}&-2.9020&&6&{\it pky $\times$ vtm}&0.6727\\
      2&{\it catct}&1.9405&&7&{\it inc $\times$ pky $\times$ vtm}&0.0801\\
      3&{\it asa}&0.3000&&8&{\it sex $\times$ inc $\times$ jomyop $\times$ asa}&0.8708\\
      4&{\it inc $\times$ pky}&0.2728&&9&{\it sex $\times$ inc $\times$ catct $\times$ asa}&0.2585\\
      5&{\it catct $\times$ asa}&0.3958\\
      \hline
    \end{tabular}
  \end{center}
  \label{mypstep1pattern}
\end{table}

\begin{figure}
  \begin{center}
    \includegraphics[scale = 0.6]{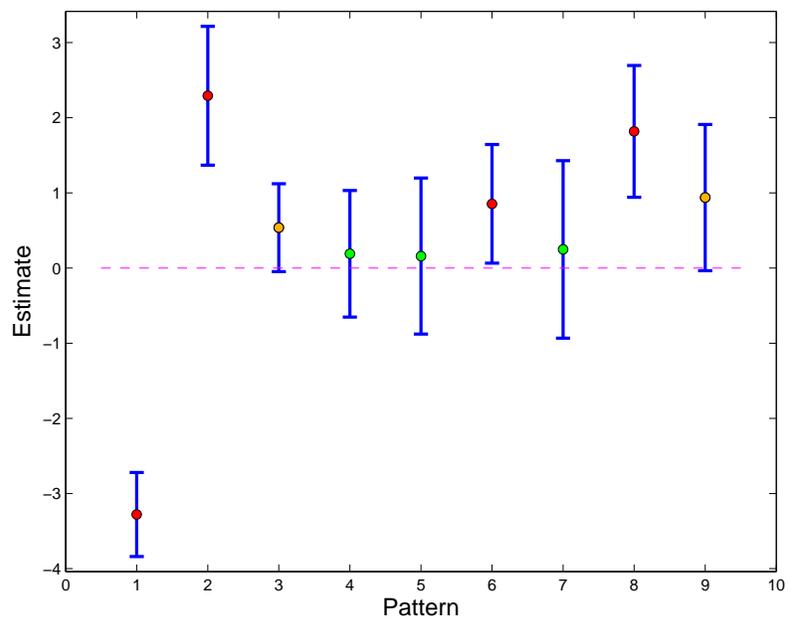}\\
  \end{center}
  \caption{
The eight patterns that survived 
Step 1 of LPS. 
The vertical bars are 90\% confidence intervals 
based on linear logistic regression.
Red dots mark the patterns that are significant
at the 90\% level. The orange dots are  borderline 
cases, the confidence intervals barely 
covering 0.}
\label{ci2}
\end{figure}

Figure \ref{ci2} plots the coefficients of the 8 patterns plus the 
constant that survived Step 1 along with 90\% confidence intervals.
These patterns are then subject to Step 2, backward elimination, tuned via 
BGACV. 
The final model after the backward elimination step is 
{{{
\begin{eqnarray}\label{myopiamodel}
f &=& -2.84 + 2.42\times catct + 1.11\times pky\times vtm \nonumber\\
  &~& + 1.98\times sex\times inc\times jomyop\times asa + 1.15\times sex\times inc\times catct\times asa. 
\end{eqnarray}
}}}
The significance levels for the coefficients of the four 
patterns in this  model (\ref{myopiamodel}) can 
be formally computed and are, respectively
3.3340e-21,
1.7253e-05,
1.5721e-04,
and 0.0428.
The pattern $pky \times vtm$ catches our attention because the pattern effect is 
strong and both variables are controllable. This model  tells us that the 
distribution of $y$, myopic change conditional on $pky = 1$ depends on 
{\it catct}, as well as {\it vtm} and higher order interactions, but 
myopic change conditional on  $pky  = 0$  is independent of {\it vtm}. This interesting effect can easily be seen by going back to a table of the original {\it catct, pky} and {\it  vtm} data (Table \ref{pkyvtmcat}). The denominators in the risk column are the 
number of persons with the given 
pattern and the numerators are the number of those with $y=1$. The first two rows list the heavy smokers with cataract. Heavy smokers  who take vitamins have a smaller risk of having myopic change. The third and fourth rows list the heavy smokers without cataract. Again, taking vitamins is protective. The first four rows suggest that taking vitamins in  heavy smokers is associated with a reduced risk of getting more myopic. The last four rows list all non-heavy smokers. Apparently taking vitamins does not similarly reduce the risk of becoming more myopic in this population. Actually, it is commonly known that smoking significantly decreases the serum and tissue vitamin level, especially Vitamin C and Vitamin E, for example  \cite{galan:viteri:bertrais:2005}. Our data suggest  a possible reduction 
in myopic change in persons who smoke who take vitamins. 
However, our data are observational and subject to 
uncontrolled confounding. A randomized controlled clinical 
trial would provide the best evidence of any effect of vitamins 
on myopic change in smokers.   
\begin{table}[htbp]
  \caption{The raw data for cataract, smoking and not taking vitamins.}
  \begin{center}
    \begin{tabular}{ccc|c}
      \hline
      catct&pky&no vitamins&risk\\\hline
      $~~~~~~~1~~~~~~~$&$~~~~~~~1~~~~~~~$&$~~~~~~~1~~~~~~~$&$~~~~~~~17/23 = 0.7391~~~~~~~$\\
      1&1&0&7/14 = 0.5000\\
      0&1&1&22/137 = 0.1606\\
      0&1&0&2/49 = 0.0408\\
      1&0&1&18/51 = 0.3529\\
      1&0&0&19/36 = 0.5278\\
      0&0&1&22/363 = 0.0606\\
      0&0&0&13/203 = 0.0640\\
      \hline
    \end{tabular}
  \end{center}
  \label{pkyvtmcat}
\end{table}

Since the model is the result of previous data mining, 
caution in making significance statements may be in 
order. To investigate the probability 
of the overall procedure to generate significant false patterns, 
we kept the attribute data fixed, randomly scrambled the 
response data and applied the LPS algorithm on the scrambled data. 
The procedure was repeated 1000 times, and in all these runs, 
1 main effect, 10 second order, 5 third order  and 
just 1 fourth order patterns showed up. 
We then checked on the raw data. There are 21 people with the pattern $sex\times inc\times jomyop\times asa$ and 9 of them have myopic change.
The incidence rate is 0.4286, as compared to the overall rate of 0.137. 
Note that none of the variables in this pattern are involved 
with variables in the two lower order patterns $catct$ and $pky\times vtm$ 
so it  can be concluded 
that the  incidence rate is contributed only by the pattern effect. 
People with the other size four pattern $sex\times inc\times catct\times asa$ 
have an incidence rate of 0.7727 (17 out of 22), which can 
be compared with the incidence rate of all people with $catct$, 0.4919.

We also applied Logic regression and SPLR on this data set. 
Logic selected both $catct$ and 
$pky\times vtm$ but missed the two high order patterns. 
Instead, it selected $asa$ as a 
main effect. Note that $asa$ is present in both size four patterns. 
SPLR selected the 
same patterns as Logic regression with an addition of $pky$, which is necessary for 
$pky\times vtm$ to be included. These results
agree with what we have found in the simulation 
studies: they are not as likely as  LPS in finding higher order patterns. 
It is noted that in the original version of LPS \cite{shi:wahba:wright:lee:2006}, 
Step 1 was tuned 
by GACV rather than BGACV, and resulted in the above eight 
patterns in Table \ref{mypstep1pattern} plus four more, 
but the final model after Step 2 was the same.
\section{Rheumatoid Arthritis and SNPs in a Generative Model Based on GAW 15}

Rheumatoid arthritis (RA) is a complex disease 
with a moderately strong genetic component. 
Generally females are at a higher risk than males.
Many studies have implicated a specific region on chromosome 6 
as being related to the risk of RA, recently \cite{thompson:worthington:2007},
although possible regions on other chromosomes have also been implicated.
The 15th Genetic Analysis Workshop 
(GAW 15, November 2006 \cite{maccluer:2007})
focused on RA, and an extensive simulation data set of 
cases and controls with simulated single nucleotide polymorphisms (SNPs)
was provided to participants and is now
publicly available \cite{miller:lind:li:jang:2007}. SNPs are DNA sequence
variations that occur when a single nucleotide in the genome 
sequence is changed.  
Many 
diseases are thought to be associated with SNP changes at multiple sites that 
may interact, thus it is important to have tools that can ferret out 
groups of possibly interacting SNPs. 

We applied LPS to some of the simulated 
SNP  RA GAW 15 data \cite{shi:lee:wahba:2007}.  
This provided an opportunity to apply LPS 
in a context with  large genetic attribute vectors, 
with a known genetic architecture, as described in 
\cite{miller:lind:li:jang:2007}, and to
compare the results 
against the description of the architecture generating the data. 
We decided to use the GAW data to build 
a simulation study where we modified the model 
that appears in \cite{shi:lee:wahba:2007}
to introduce a third order pattern, and in the 
process deal with some anomalous minus signs in
our fitted model, also observed
by others \cite{schwartz:szymczak:ziegler:konig:2007}.
We then simulated phenotype data from the GAW genotypes and 
covariates,
and can evaluate 
how well the LPS of this paper reproduces the 
model generating the data with the third order 
pattern. This section describes the results. 

In the simulated genetic data sets \cite{miller:lind:li:jang:2007} 
genome wide scans of 9187 SNPs were 
generated with just a few of the SNPs  linked to rheumatoid arthritis, 
according to the described model architecture.
The data simulation was set up to 
mimic the familial pattern of rheumatoid arthritis including a 
strong chromosome 6 effect.
A large population of nuclear families 
(two parents and two offspring) was generated. This population 
contains close to 2 million sibling pairs. From this population, a random 
sample of 1500 families was selected from among families with two 
affected offspring and another random sample of 2000 families 
was selected from among families where no member was affected. 
A total of 100 independent (replicate) data sets were generated.
We randomly picked one offspring from each family in 
replicate 1 as our data set. 
As for covariates we take the  674 SNPs in chromosome 6 that were 
generated as a subset of the genome wide scan data 
and three 
environmental variables: $age$, $sex$ and $smoking$. We created 
two dummy variables for each SNP since most of them have three 
levels: normal, one variant allele and two variant alleles. 
For environmental variables female gender is treated as a 
risk factor, smoking is treated as a risk factor and 
age $\geq 55$ is treated as a risk factor. 
We first describe a reanalysis of this data using the 
LPS algorithm of this paper. 
We began our analysis with a screen step. In this step, 
each variable is entered into a linear logistic regression model. 
For SNPs with three levels, both dummy variables are entered into the same model. 
We fit these models and keep the variables with at least one $p$-value 
less than 0.05. 
The screen step selected 74 variables (72 SNPs plus $sex$ and $smoking$). 
We then ran LPS on these 74 variables with $q = 2$, which generates 
10371 basis functions. The final model was
\begin{eqnarray}\label{gawmodel}
f& =-.62  & +~.87\times smoking + 1.05\times sex -2.04\times SNP6\_153\_1\nonumber\\
&~&-1.45\times SNP6\_154\_1+2.23\times SNP6\_162\_1 -5.60\times SNP6\_153\_2,
\end{eqnarray}
where $SNP6\_153\_1$ is SNP number 153 on chromosome 6 with 1 variant allele,  
$SNP6\_154\_1$ is SNP number 154 on chromosome 6 with 1 variant allele, $SNP6\_162\_1$ is 
SNP number 162 on chromosome 6 with 1 variant allele 
and $SNP6\_153\_2$ is SNP number 153 on chromosome 6 with 2 variant alleles. 
When the analysis in \cite{shi:lee:wahba:2007} was presented at 
the GAW 15 workshop 
in November 2006
the tuning procedure presented here had not been finalized,
and both Step 1 and 
Step 2 were tuned against prediction accuracy using replicate 2 as a 
tuning set. The model (\ref{gawmodel}) obtained here is exactly the same 
as \cite{shi:lee:wahba:2007}. We were pleased to find 
that the in-sample BGACV tuning 
here was just as good as having a separate tuning set. It is interesting
to note that in this particular problem, tuning
for prediction and for model selection apparently
led to the same results, although in general this is not necessarily 
the case. 

According to the description \cite{miller:lind:li:jang:2007}
of the architecture generating 
the data the risk of RA is affected by two loci (C and D) on chromosome 6, 
$sex$, $smoking$ and a $sex$ by locus C interaction. 
It turns out that both $SNP6\_153$ and $SNP6\_154$ are 
very close to locus C on chromosome 6 and $SNP6\_162$ is very close to locus D on chromosome 6. 
The LPS  method picked all important variables without any  false positives. 
The description of the data generation architecture said that there was 
a strong interaction between $sex$ and locus C. 
We didn't pick up a  $sex$ by locus C interaction
in \cite{shi:lee:wahba:2007}, which was surprising. 

We were curious about the apparent counter-intuitive negative 
coefficients for 
both $SNP6\_153$ patterns and the $SNP6\_154\_1$ pattern, which 
appear to say that 
normal alleles are risky and variant alleles are protective. 
Others also found an anomalous protective effect 
for $SNP6\_154$ normal alleles  \cite{schwartz:szymczak:ziegler:konig:2007}.
We went back and looked at the raw data for $SNP 6\_153$ and 
$SNP 6\_154$ as a check but 
actually Table 4 of \cite{shi:lee:wahba:2007} shows that 
this protective effect is in the simulated data, for whatever
reason, and it also shows that 
this effect is stronger for women than for men.  
We then recoded the $SNP 6\_153$ and $SNP 6\_154$ responses
to reflect the actual effect of these two variables as 
can be seen in tables of the simulated data. 

\begin{table}[htbp]
\caption{Fitted and Simulated Models, see text for explanation.
}
   \begin{center}
     \begin{tabular}{clccccc}
       \hline
       &Variable 1&Level 1&Variable 2&Level 2&Coef&Est\\\hline
       &constant&-&-&-&-4.8546&-4.6002\\
       &$smoking$&-&-&-&0.8603&0.9901\\
       Main effects&$SNP6\_153$&1&-&-&1.8911&1.5604\\
       &$SNP6\_162$&1&-&-&2.2013&1.9965\\
       &$SNP6\_154$&2&-&-&0.7700&1.0808\\\hline
       &$sex$&-&$SNP6\_153$&1&0.7848&0.9984\\
       &$sex$&-&$SNP6\_154$&2&0.9330&0.9464\\
       Second order&$SNP6\_153$&2&$SNP6\_154$&2&4.5877&4.2465\\
       patterns &$SNP6\_153$&1&$SNP6\_553$&2&0.4021&0\\
       &$SNP6\_154$&2&$SNP6\_490$&1&0.3888&0\\\hline\hline
       &             & Added             &         \\\hline
       Third order pattern&\multicolumn{4}{c}{$sex\times SNP6\_108\_2 
\times SNP6\_334\_2$}&3&2.9106\\\hline
     \end{tabular}
   \end{center}
   \label{gawmodeltab}
\end{table}

Table 6 shows the results. The new fitted model, 
above the double line, has four main effects and five 
second order patterns.
The estimated coefficients are given in the column headed
``Coef". The 
$sex \times SNP6\_153$
and $sex \times SNP6\_154$
are there as expected. Two weak second order 
patterns involving $SNP6\_553$ and $SNP6\_490$
are fitted, but do not appear to be explained by
the simulation architecture.
This model resulted from fitting with $q=2$. Then 
the LPS algorithm was run to include all third 
order patterns ($q=3$) of the 74 
variables which passed the screen step. This generated
403,594 basis functions. No third order patterns 
were found, and the fitted model was the same 
as in the $q=2$ case. To see if a third order pattern
would be found if it were there, we created a generative
model with the four main effects and five second order
patterns of Table 1,  with their coefficients
from the ``Coef" column,  
and added to it a third order pattern 
$sex\times SNP6\_108\_2 \times SNP6\_334\_2$
with coefficient 3.
The two SNPs in the third order pattern were chosen to be well 
separated in chromosome 6 from 
the reported gene loci. LPS did indeed find the 
third order pattern.The estimated coefficients
are found under the column headed ``Est". 
No noise patterns
were found, and the two weak second order patterns
in the model were missed. However, the potential
for the LPS algorithm to find higher order patterns
is clear. Further investigation of the properties 
of the method in genotype-phenotype scenarios 
is clearly warranted, taking advantage of 
the power of the LASSO algorithm to handle 
a truly large number of unknowns simultaneously. 
Run time was 4.5  minutes on an AMD Dual-Core 2.8 GHz machine
with 64 GB memory. Using multiple runs with 
clever designs to guarantee that every higher order 
pattern considered is in at least one run, will allow the analysis of 
much larger SNP data sets with tolerable computer cost. 
\section{Discussion}
In any problem where there are a large number of highly interacting predictor 
variables that are or can be reduced to dichotomous values, LPS can be 
profitably used. If the ``risky" direction 
(with respect to the outcome of interest) 
is known for all or almost 
all of the variables, the results are readily interpretable. 
If the risky direction is coded correctly for all of the 
variables, the fitted model can be expected to be 
sparser than that for any other coding. However, if a small number of 
risky variables are coded in the ``wrong" way, this usually
can be detected.
The method can be used as a preprocessor when there 
are very many continuous variables in contention, 
to reduce the number of variables for more detailed nonparametric analysis. 

LPS, using the algorithm of Appendix \ref{algorithm} is efficient.
On an Intel Quad-core 2.66 GHz machine with 8GB memory
the LPS Steps 1 and 2 can do
90,000 basis functions in 4.5 minutes. 
On an AMD Dual-Core 2.8 GHz machine with 64 GB memory
the algorithm did LPS with 403,594 basis
functions in 4.5 minutes. It  can do
2,000,000 basis functions in 1.25 hours. On the same AMD machine,
problems of the size of the myopia data (128 unknowns)
can be solved in a few seconds and the GAW 15
problem data (10371 unknowns) was solved in 1.5 minutes.

A number of considerations enter into the choice of $q$. 
If the problem is small and the user is interested in high order patterns, 
it doesn't hurt to include all possible patterns; if the problem is about 
the size of Simulation Example 3, $q=4$ might  be a good choice; 
For genomic data the choice of $q$ can be limited by extremely 
large attribute vectors.
In genomic or other data where the existence of a very small 
number of important higher patterns is suspected, but there 
are too many candidates to deal with simultaneously, it 
may be possible to overcome the curse of 
dimensionality with multiple screening levels and multiple runs. 
For example, considering say, third or even fourth 
order patterns,  variables could be assigned to doable sized runs so that
every candidate triple or quadruple of variables is assigned to at least one
run. With our purpose built algorithm, the approach is quite amenable to 
various flavors of exploratory data mining. 
When a computing system such as Condor
({\tt http://www.cs.wisc.edu/condor/})
is available, many runs can compute simultaneously.

Many generalizations are available. Two classes of models where the LASSO-Patternsearch approach 
can be expected to be useful are the multicategory end points model in \cite{linx:1998}\cite{wahba:2002b}, where an estimate is desired of the probability of being in class $k$ when there are $K$
possible outcomes; another is the multiple correlated endpoints model in \cite{gao:wahba:klein:klein:2001}. In this latter model, the correlation structure of the multiple endpoints 
can be of interest. Another generalization allows the coefficients $c_{\ell}$ to 
depend on other variables; however, the penalty functional must 
involve $\ell_1$ penalties if it is desired to have a convex optimization problem 
with good sparsity properties with respect to the patterns. 
In studies with environmental as well as genomic data 
selected interactions between SNP patterns and 
continuous covariates can be examined \cite{zhang:wahba:lin:voelker:2004}:
the numerical algorithm can be used on large collections of basis functions
that induce a reasonable design matrix, for example collections 
including splines, 
wavelets or radial basis functions.  


\section{Summary and Conclusions}

The LASSO-Patternsearch algorithm 
brings together several known ideas in a novel way, 
using a tailored tuning and pattern selection procedure and 
a new purpose built computational algorithm. We have examined the 
properties of the LPS by analysis of observational data, 
and simulation studies at a 
scale similar to the observational data. The results are verified in 
the simulation studies by examination of the generated ``truth", 
and in the observational data by selective examination of the 
observational data directly, and data scrambling to check
false alarm rates, with excellent results. The novel 
computational algorithm allows the examination of a very large 
number of patterns, and, hence, high order interactions. 
We believe the LASSO-Patternsearch will be an important addition to the toolkit of
the statistical data analyst.

\subsection*{Acknowledgments}
Thanks to David Callan for many  helpful suggestions and for Appendix D.
Research supported by 
NIH Grants EY09946 (WS, GW)
NSF Grants DMS-0505636, DMS-0604572 (WS, GW)
ONR Grant N0014-06-0095 (WS,GW)
NIH Grants  EY06594 and EY015286 (BEK, RK, KL)
Research to Prevent Blindness Senior Investigator Awards (RK, BEK)
NSF Grants SCI-0330538, DMS-0427689, CCF-0430504,
CTS-0456694, CNS-0540147, 
DOE Grant DE-FG02-04ER25627 (SW)
\appendix
\appendixpage
\addappheadtotoc

\section{The BGACV Score}\label{BGACV}
We denote the estimated logit function by 
$f_\lambda (\cdot)$ and define $f_{\lambda i} = f_\lambda(x(i))$,
$p_{\lambda i} = \frac{e^{f_{\lambda i}}}{1+e^{f_{\lambda i}}}$
$\sigma^2_{\lambda i} = 
\frac{e^{f_{\lambda i}}}{(1+e^{f_{\lambda i})^2}}$
for $i = 1,\cdots,n$. Now define
\begin{equation}
  OBS(\lambda) = \frac{1}{n}\sum_{i=1}^n[-y_if_{\lambda i} + \mbox{log}(1+e^{f_{\lambda i}})].
\end{equation}
From \cite{xiang:wahba:1996}\cite{lin:wahba:xiang:gao:klein:klein:2000}
the leave-one-out CV is
\def\fli{f_{\lambda i}}
\def\flii{f_{\lambda i}^{[-i]}}
\def\mli{p_{\lambda i}}
\def\mlii{p_{\lambda i}^{[-i]}}
\def\sn{\sum_{i=1}^{n}}
\def\si{\sigma}
\def\la{\lambda}
\begin{eqnarray} \label{cv}
CV(\lambda)
           &=&\frac{1}{n}\sum^{n}_{i=1} [-y_i f_{\la i}^{[-i]}+
                     b(f_{\la i})] \nonumber \\
	   &=& OBS(\la) + \frac{1}{n} \sn 
                [y_i(y_i -\mlii)]\left[ \frac{\fli - \flii}{y_i -\mlii}\right]\nonumber \\
           &=& OBS(\la) + \frac{1}{n}\sn 
                y_i 
		(y_i -\mli)
               \left[\frac{\fli - \flii}{y_i -\mlii}\right]
	/	
        \left[1-
        \left(\frac{\mli -\mlii}{\fli -\flii}\right)
       \left(\frac{\fli -\flii}{y_i -\mlii}\right)\right]\nonumber\\
           &\approx& OBS(\la) + \frac{1}{n}\sn 
                y_i(y_i -p_{\la i})
             \left[\frac{\fli -\flii}{y_i -\mlii}\right]/
	\left[1-\si_{\la i}^2 
		\left(\frac{\fli -\flii}{y_i -\mlii}\right)\right]\label{cvapprox}.
\end{eqnarray}
Here $\si^2_{\la i} = \si^2(f_{\la i})$ and 
the approximation in (\ref{cvapprox}) follows upon recalling that 
$\frac{\partial p}{\partial f} = \si^2$. 

Denote the objective function in (\ref{penlik}) - (\ref{J}) by $I_\lambda({y,c})$, 
let $B_{ij} =B_j(x(i))$ be the entries of the design matrix $B$,
and for ease of notation denote  
$\mu = c_{N_B}$.
Then the objective function can be written
\begin{equation}\label{penlikc}
  I_\lambda({y,c}) = \frac{1}{n}\sum_{i=1}^n[-y_i\sum_{j=1}^{N_B}c_jB_{ij} +
  {\log}(1+e^{(\sum_{j=1}^{N_B}c_jB_{ij})})] + \lambda\sum_{j=1}^{N_B -1}|c_j|.\label{cv4}
\end{equation}
\newline
\def\cl{c_{\lambda}}
\def\cle{c_{\lambda}^{\epsilon}}
\def\cls{c_{\lambda}^*)}
\def\cles{c_{\lambda}^{\epsilon *}}
Denote the minimizer of ($\mbox{\ref{cv4}}$) by $c_{\lambda}$. We know that the 
$l_1$ penalty produces sparse solutions. Without loss of generality, 
we assume that the first $s$ components of ${c_{\lambda}}$ are nonzero. 
When there is a small perturbation ${\epsilon}$ on the response, 
we denote the minimizer of $I_\lambda({c,y+\epsilon})$ by 
${c_{\lambda}^{\epsilon}}$. The 0's in the solutions are robust 
against a small perturbation in the response. That is, 
when ${\epsilon}$ is small enough, the 0 elements will stay at 0. 
This can be seen by looking at the KKT conditions when 
minimizing ($\mbox{\ref{cv4}}$). Therefore, the first 
$s$ components of $c_{\lambda}^{\epsilon}$ are nonzero and the rest are zero. 
For simplicity, we denote the first $s$ components of ${c}$ by ${c^*}$ 
and the first $s$ columns of the design matrix ${B}$ by ${B^*}$. Then
let $f_{\lambda}^y$ be the column vector with 
$i$ entry $f_{\lambda}(x(i))$ based on data $y$, and 
let $f_{\lambda}^{y + \epsilon}$ be the same column vector 
based on data $y+\epsilon$.
\begin{equation}
  {f_\lambda^{y+\epsilon} - f_\lambda^{y} = B(c_{\lambda}^{\epsilon}}-c_{\lambda})
  = B^*(c_{\lambda}^{\epsilon *}-c_{\lambda}^*).\label{cv7}
\end{equation}
Now we take the first-order Taylor expansion of 
$\frac{\partial I_{\lambda}}{{\partial c^* }}$:
\begin{equation}
  \left[\frac{\partial I_\lambda}{\partial{c^*}}\right]_{({\cle ,y+\epsilon})} 
  \approx \left[\frac{\partial I_\lambda}{{\partial c^*}}\right]_{({\cl,y})} 
  + \left[\frac{\partial^2I_\lambda}{\partial{c^*}\partial{c^{*'}}}\right]_{({\cl,y})}(\cles - \cls) 
  + \left[\frac{\partial^2I_{\lambda}}{\partial{c^*}\partial{y'}}\right]_{({\cl,y})}({y+\epsilon-y}).\label{cv5}
\end{equation}
Define
\[{U}\equiv n\left[\frac{\partial^2I_\lambda}{\partial{c^*}\partial{c^{*'}}}\right]_{({\cl,y})} = 
{{B^*}'}\mbox{diag}([\frac{e^{f_{\lambda 1}}}{(1+e^{f_{\lambda 1}})^2},\cdots,\frac{e^{f_{\lambda n}}}{(1+e^{f_{\lambda n}})^2}]){B^*} 
= B*'WB*, ~~\mbox{say}, \]
and
\[{V} \equiv -n\left[\frac{\partial^2I_{\lambda}}{\partial{c^*}\partial{y'}}\right]_{({\cl,y})} = {{B^*}'}.\]
By the first-order conditions, the left-hand side 
and the first term of the right-hand side of ($\mbox{\ref{cv5}}$) are zero. So we have
\begin{equation}
U(\cles - \cls)\approx {V\epsilon}.\label{cv6}
\end{equation}
Combine ($\mbox{\ref{cv7}}$) and ($\mbox{\ref{cv6}}$) we have ${f_\lambda^{y+\epsilon}-f_\lambda^y} \approx {H\epsilon}$, where 
\begin{equation}\label{H}
{H \equiv B^*U^{-1}V} \equiv B^*U^{-1}{B^*}'.
\end{equation} 
Now let $\epsilon$ be    
${\epsilon}=(0,\cdots,y_i - p_{\lambda i}^{[-i]},\cdots,0)'$; 
then ${f_{\lambda}^{y+\epsilon}-f_\lambda^y}\approx{H}_{i}\epsilon_i$, 
where $\epsilon_{i} = y_i -p_{\lambda i}^{[-i]}$ 
and $H_{i}$ is the $i$th column of $H$. 
By the Leave-Out-One Lemma (stated below), 
${f_\lambda}^{[-i]}$ = ${f_\lambda}^{{y+\epsilon}}$. Therefore
\begin{equation}
   \frac{f_{\lambda i}-f_{\lambda i}^{[-i]}}{y_i - p_{\lambda i}^{[-i]}}
  = \frac{f_{\lambda i}^{y+\epsilon}-f_{\lambda i}}{y_i -p_{\lambda i}^{[-i]}}\approx h_{ii}\nonumber
\end{equation}
where $h_{ii}$ is the $ii$th entry of $H$. 
From the right hand side of ($\mbox{\ref{cvapprox}}$), 
the approximate CV score is
\begin{equation}\label{ACV}
  CV(\lambda) \approx \frac{1}{n}\sum_{i=1}^n[-y_if_{\lambda i} + 
\mbox{log}(1+e^{f_{\lambda i}})] 
+ \frac{1}{n}\sum_{i=1}^nh_{ii}
\frac{y_i(y_i-p_{\lambda i})}{(1-\sigma_{\lambda i}^2h_{ii})}\label{acv}.
\end{equation}
The $GACV$ score is obtained from the approximate CV score in (\ref{ACV}) by replacing 
$h_{ii}$ by $\frac{1}{n} tr(H)$ and $\sigma_{\lambda i}^2 h_{ii}$ by $\frac{1}{n} tr(WH)$.
It is not hard to see that $tr(WH) = trW^{1/2}H W^{1/2} = s \equiv N_{B_0}$,
the number of basis functions 
in the model, giving 

\begin{equation}\label{gacvA}
GACV(\lambda) = \frac{1}{n}\sum_{i=1}^n[-y_if_{\lambda i} 
+ \mbox{log}(1+e^{f_{\lambda i}})] + 
\frac{1}{n}
{\mbox{tr}}H
\frac{\sum_{i=1}^ny_i(y_i-p_{\lambda i})}{(n - N_{B_0}) }.
\end{equation}
Adding the weight $\frac{1}{2} \log n $ to 
the ``optimism'' part of the GACV score, we obtain the 
B-type GACV (BGACV):
\begin{equation}\label{bgacvA}
BGACV(\lambda) = \frac{1}{n}\sum_{i=1}^n[-y_if_{\lambda i} 
+ \mbox{log}(1+e^{f_{\lambda i}})] + 
\frac{1}{n}
\frac{\log n}{2}\mbox{tr}H
\frac{\sum_{i=1}^ny_i(y_i-p_{\lambda i})}{(n - N_{B_0}) }.
\end{equation}

\begin{lemma}\label{lem1}{(Leave-Out-One Lemma)}\\
Let the objective function $I_{\lambda}(y,f)$ be defined as 
before.
Let $f_{\lambda}^{[-i]}$ be the minimizer of $I_\lambda(y,f)$ with the $i$th 
observation omitted and let $p_\lambda^{[-i]}$ be the corresponding probability. 
For any real number  $\nu$, we define the vector 
$z = (y_1,\cdots,y_{i-1},\nu,y_{i+1},\cdots,y_n)'$. 
Let $h_\lambda(i,\nu,\cdot)$ be the minimizer of $I_\lambda({z, f})$; then 
$h_\lambda(i,p_\lambda^{[-i]},\cdot)=f_\lambda^{[-i]}(\cdot)$.
\end{lemma}
The proof of lemma $\ref{lem1}$ is quite simple and very similar to the proof of 
the Leave-Out-One-Lemma in \cite{zhang:wahba:lin:voelker:2004} so we will omit it here. 

We remark that in this paper we have employed the BGACV criterion
twice as a stringent model selector under the assumption that
the true or the desired model is sparse. Simulation experiments
(not shown) suggest that the GACV criterion is preferable if the true model
is not sparse and/or the signal is weak. 
The GACV and the BGACV selections  probably 
bracket the region of interest of $\lambda$ in most applications. 
\subsection*{Comments}
A referee has asked how BGACV might be compared to 
the more familiar 
$BIC = - \log \textrm{likelihood} + \frac{\log n}{2} df$.
where $df$ is the degrees of freedom in the case of 
Bernoulli data. The short answer to this question 
is that an exact expression 
for $df$ does not, in the usual sense, exist in the 
case of Bernoulli data. Thus, only a hopefully
good approximation to something
that plays the role of $df$ in the 
Bernoulli case can be found. 
This argument, which is independent of the nature 
of the estimate $f_{\lambda}$  of $f$,
is found in Section 2 of 
\cite{lin:wahba:xiang:gao:klein:klein:2000}.
We  sketch the main idea.
Let $KL(\lambda) = KL(f,f_{\lambda})$ be the 
Kullback-Liebler distance between the distribution
with the true but unknown canonical link $f$ and the 
distribution with link 
$f_{\la}$ and let  
\begin{equation}
CKL(\la) =  KL(f,f_{\la}) -\frac{1}{n}[
\sum_{i=1}^n E_{\mu}y_i f_i -b(f_i)] 
        \equiv \sum_{i=1}^n -\mu_i f_{\la i} + b(f_{\la i})\label{ckl}
\end{equation}
be the comparative Kullback-Liebler distance.
The goal is to find 
an unbiased estimate of $CKL(\lambda)$ as a function 
of $\lambda$, which will then be minimized to estimate the 
$\lambda$ minimizing the true but unknown  $CKL$. Letting 
$OBS(\la) = \frac{1}{n} \sum_{i=1}^n -y_i f_{\lambda i} + b(f_{\lambda i})$
we can write 
$CKL(\lambda)  = OBS(\lambda) + D(\lambda)$.
where
$D(\lambda) = \frac{1}{n}\sum_{i=1}^n(y_i -\mu_i) f_{\lambda i}$.
Then 
$E_{\mu}D(\la) = \frac{1}{n} \sum_{i=1}^n E_{\mu}(y_i -\mu_i)f_{\lambda i}$.
Ye and Wong (1997b) show, in exponential families,
for {\it any} estimate $f_{\lambda}$ of $f$
\begin{equation}\label{EU}
\frac{1}{n}\sum_{i=1}^n E_{\mu_i}(y_i -\mu_i) f_{\lambda i}
=\frac{1}{n}\sum_{i=1}^n \sigma^2_i
\frac{\partial}{\partial \mu_i} E_{\mu_i}({f_{\lambda i}}).
\end{equation}
Here $E_{\mu_i}(f_{\lambda i})$ is the expectation with 
respect to $y_i$ conditional on the  $y_j, j \neq i$ being
fixed. 
(Their proof is reproduced in \cite{lin:wahba:xiang:gao:klein:klein:2000}.)
Ye and Wong call $n$ times the right hand side of 
(\ref{EU}) the generalized degrees of freedom (GDF). 
and it does indeed reduce
to the usual trace of the influence matrix in the case of 
Gaussian data with quadratic penalties. 
Unbiased estimates of the GDF can be found for 
Poisson, Gamma, Binomial distribution taking on three or more values, 
and other distributions, using the results in \cite{hudson:1978} 
but Ye and Wong show that {\it no unbiased estimate of the GDF 
in the Bernoulli case exists.} Thus, in the absence
of a bona fide unbiased risk method of estimating the 
$CKL(\lambda)$ the alternative GACV based on leave-one-out 
to target the $CKL$ has been proposed. Thus, 
in this paper 
$\frac{1}{n} \mbox{tr}H
\frac{\sum_{i=1}^ny_i(y_i-p_{s i})}{(n - N_{B_0}) }
= \hat{D}(\lambda)$, say, plays the role of $df$. 
Since no exact unbiased estimate for  
$df$ exists, the issue of the accuracy of the approximations
in obtaining $\hat{D}$ reduces to the issue of to what extent
the minimizer of $GACV(\lambda)$ is a good estimate of 
the minimizer of the (unobservable) $CKL(\lambda)$. 
The GACV for Bernoulli data was first proposed in 
\cite{xiang:wahba:1996} for RKHS (quadratic) penalty functionals, 
where simulation results demonstrated the accuracy of this
approximation. Further excellent favorable results 
for a randomized version of the GACV with RKHS penalties were presented in
\cite{lin:wahba:xiang:gao:klein:klein:2000}. 
In \cite{zhang:wahba:lin:voelker:2004} the GACV was derived 
for Bernoulli data with  $l_1$ penalties with a nontrivial 
null space, and favorable results for the randomized 
version were obtained. 
The derivation was rather complicated, and a simplified
derivation as well as a simpler expression for the result which
is possible in the present context  
are presented above. 
A recent
work involving the LASSO in the Bernoulli case with $l_1$ penalty
uses a tuning set to choose the smoothing parameters. 
SPLR \cite{park:hastie:2008} uses 
$tr(B*'WB* + \lambda I*)^{-1}(B*'WB*)$,
where $I$ is the diagonal matrix with 
all 1's except in the position of the model 
constant, 
as their proxy for df in their BIC-like criteria
for model selection after fixing $\lambda$.
In the light of the Ye and Wong result it is no 
surprise that an exact definition of $df$ in this case 
cannot be found in the literature.
\section{Minimizing the Penalized Log Likelihood Function}
\label{algorithm}
The function (\ref{penlik}) is not differentiable with respect to the
coefficients $\{c_{\ell}\}$ in the expansion (\ref{J}), so most
software for large-scale continuous optimization cannot be used to
minimize it directly. We can however design a specialized algorithm
that uses gradient information for the smooth term ${\cal L}(y,f)$ to
form an estimate of the correct active set (that is, the set of
components $c_{\ell}$ that are zero at the minimizer of
(\ref{penlik})). Some iterations of the algorithm also attempt a
Newton-like enhancement to the search direction, computed using the
projection of the Hessian of ${\cal L}$ onto the set of nonzero
components $c_{\ell}$. This approach is similar to the two-metric
gradient projection approach for bound-constrained minimization, but
avoids duplication of variables and allows certain other economies in
the implementation.

We give details of our approach by simplifying the notation and
expressing the problem as follows:
\begin{equation} \label{min.T}
\min_{z \in R^m} \, T_{\lambda}(z) := T(z) + \lambda \| z \|_1.
\end{equation}
When $T$ is convex (as in our application), $z$ is optimal for
(\ref{min.T}) if and only if the following condition holds:
\begin{equation} \label{opt.T}
\nabla T(z) + \lambda v = 0,
\end{equation}
for some vector  $v$ in the subdifferential of $\| z\|_1$ (denoted by
$\partial \| z\|_1$), that is,
\begin{equation} \label{subdiff.z1}
v_i \left\{ \begin{array}{ll}
 = -1 & \;\; \makebox{\rm if $z_i<0$} \\
 \in [-1,1] & \;\; \makebox{\rm if $z_i=0$} \\
 = 1 & \;\; \makebox{\rm if $z_i>0$} 
\end{array} \right.
\end{equation}
A measure of near-optimality is  given as follows:
\begin{equation} \label{nearopt}
\delta(z) = \min_{v \in \partial \|z\|_1} \| \nabla T(z) + \lambda v \|.
\end{equation}
We have that $\delta(z)=0$ if and only if $z$ is optimal.

In the remainder of this section, we describe a simplified version of
the algorithm used to solve (\ref{min.T}), finishing with an outline
of the enhancements that were used to decrease its run time.

The basic (first-order) step at iteration $k$ is obtained by forming a
simple model of the objective by expanding around the current iterate
$z^k$ as follows:
\begin{equation} \label{T.1os}
d^k = \arg \min_d \, T(z^k) + \nabla T(z^k)^Td + \frac12 \alpha_k d^Td + 
\lambda  \| z^k+d \|_1,
\end{equation}
where $\alpha_k$ is a positive scalar (whose value is discussed below)
and $d^k$ is the proposed step.  The subproblem (\ref{T.1os}) is
separable in the components of $d$ and therefore trivial to solve in
closed form, in $O(m)$ operations.  We can examine the solution $d^k$
to obtain an estimate of the active set as follows:
\begin{equation} \label{actset.est}
{\cal A}_k = \{ i = 1,2,\dots,m \; | \; (z^k+d^k)_i = 0 \}.
\end{equation}
We define the ``inactive set'' estimate ${\cal I}_k$ to be the
complement of the active set estimate, that is,
\[
{\cal I}_k = \{ 1,2,\dots,m \} \setminus {\cal A}_k.
\]

If the step $d^k$ computed from (\ref{T.1os}) does not yield a
decrease in the objective function $T_{\lambda}$, we can increase
$\alpha_k$ and re-solve (\ref{T.1os}) to obtain a new $d^k$. This
process can be repeated as needed. It can be shown that, provided
$z^k$ does not satisfy an optimality condition, the $d^k$ obtained
from (\ref{T.1os}) will yield $T_{\lambda}(z^k+d^k)<T_{\lambda}(z^k)$
for $\alpha_k$ sufficiently large.

We enhance the step by computing the restriction of the Hessian
$\nabla^2 T(z^k)$ to the set ${\cal I}_k$ (denoted by
$\nabla^2_{{\cal I}_k {\cal I}_k} T(z^k)$) and then computing a
Newton-like step in the ${\cal I}_k$ components as follows:
\begin{equation} \label{red.Newton}
( \nabla^2_{{\cal I}_k {\cal I}_k} T(z^k) + \delta_k I) p^k_{{\cal I}_k} = 
- \nabla_{{\cal I}_k} T(z^k) - \lambda w_{{\cal I}_k},
\end{equation}
where $\delta_k$ is a small damping parameter that that goes to zero
as $z^k$ approaches the solution, and $w_{{\cal I}_k}$ captures the
gradient of the term $\|z\|_1$ at the nonzero components of $z^k+d^k$.
Specifically, $w_{{\cal I}_k}$ coincides with $\partial \|z^k+d^k
\|_1$ on the components $i \in {\cal I}_k$.  If $\delta_k$ were set to
zero, $p^k_{{\cal I}_k}$ would be the (exact) Newton step for the
subspace defined by ${\cal I}_k$; the use of a damping parameter
ensures that the step is well defined even when the partial Hessian
$\nabla^2_{{\cal I}_k {\cal I}_k} T(z^k)$ is singular or nearly
singular, as happens with our problems.  In our implementation, we
choose
\begin{equation} \label{damping}
\delta_k = \min \left( \delta(z^k), \makebox{\rm mean diagonal of } 
\nabla^2_{{\cal I}_k {\cal I}_k} T(z^k) \right),
\end{equation}
where $\delta(z)$ is defined in  (\ref{nearopt}).

Because of the special form of $T(z)$ in our case (it is the function
${\cal L}$ defined by (\ref{defL}) and (\ref{f.def})), the Hessian is
not expensive to compute once the gradient is known. However, it is
dense in general, so considerable savings can be made by evaluating
and factoring this matrix on only a reduced subset of the variables,
as we do in the scheme described above.

If the partial Newton step calculated above fails to produce a
decrease in the objective function $T_{\lambda}$, we reduce its length
by a factor $\gamma_k$, to the point where $z^k_i+\gamma_k p^k_i$ has
the same sign as $z^k_i$ for all $i \in {\cal I}_k$. If this modified
step also fails to decrease the objective $T_{\lambda}$, we try the
first-order step calculated from (\ref{T.1os}), and take this step if
it decreases $T_{\lambda}$. Otherwise, we increase the parameter
$\alpha_k$, leave $z^k$ unchanged, and proceed to the next iteration.

We summarize the algorithm as follows.

\begin{algorithm}$~$\\
{\bf given} initial point $z^0$, initial damping $\alpha_0>0$,
constants ${\tt tol}>0$ and $\eta \in (0,1)$; \\
{\bf for} $k=0,1,2, \dots$ \\
\indent {\bf if} $\delta(z^k)< {\tt tol}$ \\
\indent \indent {\bf stop} with approximate solution $z^k$; \\
\indent {\bf end} \\
\indent Solve (\ref{T.1os}) for $d^k$;  (* first-order step *) \\
\indent Evaluate ${\cal A}_k$ and ${\cal I}_k$; \\
 \indent Compute $p^k_{{\cal I}_k}$ from (\ref{red.Newton}); 
(* reduced Newton step *) \\
 \indent Set $z^+_{{\cal I}_k} = z^k_{{\cal I}_k} + p^k_{{\cal I}_k}$ and
$z^+_{{\cal A}_k}=0$; \\
 \indent {\bf if} $T_{\lambda}(z^+) < \min (T_{\lambda}(z^k+d^k),T_{\lambda}(z^k))$ (* Newton step successful *) \\
 \indent \indent $z^{k+1} \leftarrow z^+$; \\
 \indent {\bf else} \\
 \indent \indent Choose $\gamma_k$ as the largest positive number such that $(z^k+\gamma_k p^k)_i  z^k_i>0$ \\
 \indent \indent \indent for all $i$ with $z^k_i \neq 0$; 
(* damp the Newton step *) \\
 \indent \indent Set $z^+_{{\cal I}_k} = z^k_{{\cal I}_k} + \gamma_k p^k_{{\cal I}_k}$ and $z^+_{{\cal A}_k}=0$; \\
 \indent \indent {\bf if}  $T_{\lambda}(z^+) < \min (T_{\lambda}(z^k+d^k),T_{\lambda}(z^k))$ (* damped Newton step successful *) \\
 \indent \indent \indent $z^{k+1} \leftarrow z^+$; \\
 \indent \indent {\bf else} {\bf if} $T_{\lambda}(z^k+d^k) < T_{\lambda}(z^k)$ 
(* first-order step successful; use it if Newton steps have failed *) \\
 \indent \indent \indent $z^{k+1} \leftarrow z^k + d^k$; \\
 \indent \indent {\bf else} (* unable to find a successful step *) \\
 \indent \indent \indent $z^{k+1} \leftarrow z^k$; \\
 \indent \indent {\bf end} \\
 \indent {\bf end} \\
 \indent (* increase or decrease $\alpha$ depending on success of first-order step *) \\
 \indent {\bf if} $T_{\lambda}(z^k+d^k) < T_{\lambda}(z^k)$ \\
 \indent \indent $\alpha_{k+1} \leftarrow \eta \alpha_k$; 
(* first-order step decreased $T_{\lambda}$, so decrease $\alpha$ *) \\
 \indent {\bf else}  \\
 \indent \indent $\alpha_{k+1} \leftarrow \alpha_k / \eta$; \\
 \indent {\bf end} \\
{\bf end}
\end{algorithm}

We conclude by discussing some enhancements to this basic approach
that can result in significant improvements to the execution time.
Note first that evaluation of the full gradient $\nabla T(z^k)$, which
is needed to compute the first-order step (\ref{T.1os}) can be quite
expensive. Since in most cases the vast majority of components of
$z^k$ are zero, and will remain so after the next step is taken, we
can economize by selecting just a subset of components of $\nabla
T(z^k)$ to evaluate at each step, and allowing just these components
of the first-order step $d$ to be nonzero. Specifically, for some
chosen constant $\sigma \in (0,1]$, we select $\sigma m$ components
  from the index set $\{ 1,2,\dots,m\}$ at random (using a different
  random selection at each iteration), and define the working set
  ${\cal W}_k$ to be the union of this set with the set of indices $i$
  for which $z^k_i \neq 0$. We then evaluate just the components of
  $\nabla T(z^k)$ for the indices $i \in {\cal W}_k$, and solve
  (\ref{T.1os}) subject to the constraint that $d_i=0$ for $i \notin
  {\cal W}_k$.

Since $\delta(z^k)$ cannot be calculated without knowledge of the full
gradient $\nabla T(z^k)$, we define a modified version of this
quantity by taking the norm in (\ref{nearopt}) over the vector defined
by ${\cal W}_k$, and use this version to compute the damping parameter
$\delta_k$ in (\ref{damping}).

We modify the convergence criterion by forcing the {\em full} gradient
vector to be computed on the next iteration $k+1$ when the threshold
condition $\delta(z^k)< {\tt tol}$ is satisfied. If this condition is
satisfied again at iteration $k+1$, we declare success and terminate.

A further enhancement is that we compute the second-order enhancement
only when the number of components in ${\cal I}_k$ is small enough to
make computation and factorization of the reduced Hessian
economical. In the experiments reported here, we compute only the
first order step if the number of components in ${\cal I}_k$ exceeds
$500$.
\newpage

\section{Results of Simulation Example 3}
\label{appsim3}
\begin{table}[htbp]
 \caption{Results of Simulation Example 3. In each row of a cell the first three numbers are the appearance frequencies of the important patterns and the last number is the appearance frequency of noise patterns.}
 \begin{center}
    \begin{tabular}{l|c|c|c|c|c}
      \hline
      $\rho_2\backslash\rho_1$&&{0}&{0.2}&{0.5}&{0.7}\\
      \hline
      &LPS&96/100/100/54&98/100/100/46&100/100/100/43&100/100/100/44\\
      0&Logic&100/98/96/120&98/95/93/107&99/94/92/83&100/98/83/134\\
      &SPLR&100/100/100/527&100/100/98/525&100/100/98/487&100/100/97/489\\\hline
      &LPS&99/100/100/46&100/100/100/49&100/100/100/39&100/100/98/36\\
      0.2&Logic&99/97/94/96&100/99/87/94&100/100/88/73&100/99/86/117\\
      &SPLR&100/100/94/517&100/99/96/530&100/97/95/495&100/100/96/485\\\hline
      &LPS&99/100/99/47&99/100/100/51&100/100/99/51&100/100/98/46\\
      0.5&Logic&99/96/86/162&100/95/87/109&100/96/78/122&100/99/80/143\\
      &SPLR&100/98/75/548&100/96/80/552&100/99/80/531&100/98/78/518\\\hline
      &LPS&100/99/96/44&99/99/97/51&100/99/96/67&100/99/94/65\\
      0.7&Logic&100/83/70/195&100/88/69/167&100/85/70/153&100/89/74/126\\
      &SPLR&100/91/51/580&100/85/49/594&100/81/52/584&100/72/55/570\\\hline
    \end{tabular}
  \end{center}
    \label{sim3tab}
\end{table}

\section{Effect of Coding Flips}\label{codeflip}
Proposition: Let $f(x) = \mu + \sum c_{j_1j_2..j_r}B_{j_1j_2..j_r}(x)$
with all $c_{j_1j_2..j_r}$ which appear in the 
sum strictly positive. If $x_j \rightarrow 1-x_j$ for 
$j \in$ some subset of $\{1,2, ..., p\}$ such that at least one 
$x_j$ appears in $f$, then the resulting representation has 
at least one negative coefficient and at least as many terms as $f$.
This follows from the \\

Lemma: Let $g_k(x)$ be the function obtained from 
$f$ by transforming $x_j \rightarrow 1-x_j, 1 \leq j \leq k$. 
Then the coefficient of $B_{j_1j_2..j_r}(x)$ in 
$g_k(x)$ is 
\[ (-1) ^{\vert \{j_1, .., j_r\} \cap \{ 1, ..., k \}\vert}
\sum_{ \{j_1, \cdots, j_r \} \subseteq T \subseteq \{j_1, \cdots, j_r\}\cup \{1, \cdots k,\}}
c_T
\]
where $\vert \cdot\vert$ means number of entries. 
\bibliographystyle{plain}

\end{document}